\documentclass[a4paper, 12pt]{amsart}
\numberwithin{equation}{section}
\usepackage{stackrel}
\usepackage{amssymb}
\usepackage{amsmath}
\usepackage{amsfonts}
\usepackage{marvosym}
\usepackage{tikz}
\newcommand*\circled[1]{\tikz[baseline=(char.base)]
{\node[shape=circle,draw,inner sep=1] (char) {#1};}}
\usepackage{stackrel}

\newtheorem{theorem}{Theorem}[section]
\newtheorem{thm}{Theorem}[section]
\newtheorem{lem}[theorem]{Lemma}

\newtheorem{prop}[theorem]{Proposition}

\theoremstyle{definition}

\newtheorem{rem}[theorem]{Remark}

\def\ca{{\mathcal A}}
\def\cb{{\mathcal B}}

\def\ce{{\mathcal E}}

\def\ch{{\mathcal H}}

\def\ck{{\mathcal K}}

\def\cam{{\mathcal M}}
\def\cn{{\mathcal N}}

\def\cs{{\mathcal S}}

\def\ga{{\mathfrak A}}
\def\gb{{\mathfrak B}}

\def\gam{{\mathfrak M}}

\def\bc{{\mathbb C}}

\def\bm{{\mathbb M}}
\def\bn{{\mathbb N}}
\def\bp{{\mathbb P}}

\def\bz{{\mathbb Z}}

\def\a{\alpha}
\def\b{\beta}
\def\g{\gamma}  \def\G{\Gamma}
  
\def\eeps{\epsilon}

\def\l{\lambda} 
\def\k{\kappa}
\def\m{\mu}

\def\n{\nu}
\def\r{\rho}
\def\s{\sigma} 

\def\f{\varphi}  
\def\th{\theta} 
\def\om{\omega}

\def\ots{\overline{\otimes}}

\def\tr{\mathop{\rm Tr}}

\newcommand{\ty}[1]{\mathop{\rm {#1}}}
\def\di{{\rm d}}

\def\ad{\mathop{\rm ad}}

\def\idd{{1}\!\!{\rm I}}

\DeclareMathAlphabet{\mathpzc}{OT1}{pzc}{m}{it}

\begin{document}

\title[klein transformation]
{Symmetric states for $C^*$-Fermi systems II: Klein transformation and their structure}
\author{Francesco Fidaleo}
\address{Francesco Fidaleo\\
Dipartimento di Matematica \\
Universit\`{a} di Roma Tor Vergata\\
Via della Ricerca Scientifica 1, Roma 00133, Italy} \email{{\tt
fidaleo@mat.uniroma2.it}}
\date{\today}

\begin{abstract}
In the present note, which is the second part of a work concerning the study of the set of the symmetric states,
we introduce the extension of the Klein transformation for general Fermi tensor product of two $\bz^2$ graded $C^*$-algebras, under the condition that the grading of one of the involved algebras is inner. After extending the construction to $C^*$-inductive limits, such a Klein transformation 
realises a canonical $*$-isomorphism between two $\bz^2$-graded $C^*$-algebras made of
the infinite Fermi $C^*$-tensor product 
$$
\ga_{\rm F}:=\big(\circled{\rm{{\tiny F}}}_{\,\bn}\, \gb, \circled{\rm{{\tiny F}}}_{\,\bn}\, \a\big)\,,
$$ 
and the infinite $C^*$-tensor product 
$$
\ga_{\rm X}:=\big(\circled{\rm{{\tiny X}}}_{\,\bn}\, \gb, \circled{\rm{{\tiny X}}}_{\,\bn}\, \a\big)
$$
of a single $\bz^2$-graded $C^*$-algebra $(\gb,\a)$, both built with respect to the corresponding minimal $C^*$-cross norms.
It preserves the grading, and its transpose sends even product states of $\ga_{\rm X}$ in (necessarily even) product states on $\ga_{\rm F}$, and therefore induces an isomorphism of simplexes
$$
\cs_\bp(\ga_{\rm F})=\cs_{\bp\times\bz^2}(\ga_{\rm F})\sim\cs_{\bp\times\bz^2}(\ga_{\rm X})\,,
$$
which allows to reduce the study of the structure of the symmetric states for $C^*$-Fermi systems to the corresponding even symmetric states on the usual infinite $C^*$-tensor product. Other relevant properties of symmetric states on the Fermi algebra will be proved without the use of the Klein transformation.

We end with an example for which such a Klein transformation is not implementable, simply because the Fermi tensor product does not generate a usual tensor product. Therefore, in general, the study of the symmetric states on the Fermi algebra cannot be reduced to that of the corresponding symmetric states on the usual infinite tensor product, even if both share many common properties.

\end{abstract}

\maketitle

\section{Introduction}
\label{sec1}

The systematic investigation of systems describing Fermi particles or, in other words, having a natural $\bz^2$ symmetry describing the {\it univalence superselection rule} ({\it e.g.} \cite{SW}), if was natural for the applications to quantum field theories, on the other side it
has had an impetuous growth in the last decades for the applications to other fields of the mathematical physics, like spin-like Fermi models on the lattices, see {\it e.g.} \cite{NSY} and the references cited therein. The study of disordered Fermi systems, and quantum Markov stochastic processes enjoying the $\bz^2$-symmetry, was recently also addressed in \cite{BF, AFM, F1}.\footnote{A kind of Lieb-Robinson bound for Fermi models on lattices was firstly established in \cite{BF}, Theorem 7.4.} 

All these various and interesting applications of Fermi systems, for which the combinatorics necessary to handle the $\bz^2$ symmetry, introduces not simply difficulties to overcome, suggested to address the systematic study  $\bz^2$-graded $C^*$-algebras in a purely abstract way. 

This project started in \cite{CDF} with the aim to extend and investigate the {\it quantum detailed balance} for $\bz^2$ graded systems, both from a general point of view and for applications to concrete models involving the CAR algebra. The key-point for this investigation was a kind of the so-called Klein transformation, which was previously applied essentially for concrete $\bz^2$-graded systems represented on some Hilbert space. To be more precise, a {\it Klein transformation} is a $*$-isomorphism which realises a similarity between the Fermi and the usual tensor product systems.

The Klein transformation was a very powerful tool for the investigations of models in quantum field theories, again all faithfully represented on some Hilbert space but, except the simple case of the usual CAR algebra, no general description of such a fundamental transformation was available up to now.

The present note is the second part relative to the investigation of the main properties of the so called symmetric states on Fermi systems. 

Indeed, to be more precise, in the 1st part \cite{Fde}, 
we built the infinite $C^*$-Fermi tensor product 
$(\ga_{\rm F},\a):=\big(\circled{\rm{{\tiny F}}}_{\,\bn}\, \gb, \circled{\rm{{\tiny F}}}_{\,\bn}\, \b\big)$ of the single $\bz^2$-graded $C^*$-algebra $\gb$ by using the analogous of the minimal norm of the usual $C^*$-tensor product firstly introduced at the end of Section 9 of \cite{CDF}.

Such a Fermi algebra $\ga_{\rm F}$ is equipped with a $\bz^2$-grading obtained by the direct limit of the grading of $\b$ of $\gb$. Notice that, the usual infinite $C^*$-tensor product is also a $\bz^2$-graded in a natural way, which is induced by the infinite tensor product of the grading $\b$ of
$\gb$. 

As for the simplest case of the CAR algebra, in the present paper we provide the definition of the Klein transformation in this general setting, under the additional condition of the innerness of the grading automorphism. After a painstaking analysis, we show that such a Klein transformation realises a $*$-isomorphism between $\ga_{\rm F}$ and $\ga_{\rm X}$ preserving the grading, and whose transposed map sends product states on a single even one on $\gb$ for $\ga_{\rm X}$, to the analogous one on $\ga_{\rm F}$.

As an almost immediate consequence, the transposed map of the Klein transformation realises an equivalence 
$$
\cs_\bp(\ga_{\rm F})=\cs_{\bp\times\bz^2}(\ga_{\rm F})\sim\cs_{\bp\times\bz^2}(\ga_{\rm X})
$$
(here $\bp$ denotes the group of all finite permutations of the natural numbers, and $\bz^2$ takes into account of the parity)
of simplexes which allows us to investigate the structure of the symmetric states on the Fermi algebra $\ga_{\rm F}$ in terms of the corresponding even and symmetric ones on the usual infinite tensor product $\ga_{\rm X}$. This is in fact the last part of the present work in which many results contained in the seminal paper \cite{St2}, and relative to $\ga_{\rm X}$, can be rephrased for $\ga_{\rm F}$ by using such an identification.

We also prove most on the remaining main results relative to symmetric states on the Fermi algebra $\ga_{\rm F}$ by using the ergodic properties of the $C^*$-dynamical system $(\ga_{\rm F},\bp)$ shared with the analogous ones of  
$(\ga_{\rm F},\bp\times\bz^2)$. In so doing, we restore a similarity between $\cs_\bp(\ga_{\rm F})$ and $\cs_{\bp\times\bz^2}(\ga_{\rm X})$ involving many basic properties, even when the Klein transformation is not implementable.

We end the introduction by remarking the following facts. Due to the extremely high degree of entanglement introduced by the (anti)com\-mutation relations, the
involved GNS representations do not generate a usual tensor product, even if the involved states are product states. This means that a Klein transformation cannot be implementable in general.

Nevertheless, we would like to remark the interest in itself of the abstract version of the Klein transformation introduced an studied here, and its potential applications to quantum field theory and other fields like
$\bz^2$-symmetric models of statistical mechanics and quantum probability.

Concerning the comparison of symmetric states on the Fermi algebra and the corresponding one on the usual infinite tensor product, we show with a simple but illustrative example that the Klein transformation cannot indeed be implemented in the full generality.

\section{Preliminaries}

We gather here some facts useful for the forthcoming sections, and refer the reader to Section 2 of \cite{Fde} for the remaining notations and basic results frequently used here. 

In the present paper all involved $C^*$-algebras are unital, equipped with a $\bz^2$-grading, and the used norm on the $C^*$-Fermi tensor products  and the usual ones are always the minimal $C^*$-cross norms if it is not otherwise specified. No assumption about the separability of the involved $C^*$-algebras is done if it is not explicitly specified.
\medskip

\noindent
\textbf{Basic notions.} 
Let $\ga$ be a $C^*$-algebra, and $\f\in\cs(\ga)$. Consider the support-projection of the canonical image of $\f$ in the enveloping von Neumann algebra $\ga^{**}$ ({\it e.g.} \cite{T}, Section III.2). It is denoted by $s(\f)$ with an abuse of notation.
\begin{prop}[\cite{NSZ}, pag. 15]
The state $\f\in\cs(\ga)$ has central support if and only if the cyclic vector $\xi_\f$ of the GNS representation of $\f$ is also cyclic for $\pi_\f(\ga)'$.
\end{prop}

For a $\bz^2$-graded $C^*$-algebra $(\ga,\a)$, we denote by 
$$
\pi_\a:=\bigoplus_{\f\in\cs_+(\ga)}\pi_\f\
$$
its universal representation built by using the even states w.r.t. the grading induced by $\a$, and thus "$\a$" stands only to recalling that $\pi_\a$ is built taking into account which grading we are using.

For the convenience of the reader, we also report the following well-known result asserting that such representations $\pi_\a$ are indeed universal 
representations for $\bz^2$-graded algebras $\ga$, independently on the chosen grading automorphism $\a$.
\begin{lem}
\label{unpasd}
The even states of a $\bz^2$-graded $C^*$-algebra $(\ga,\a)$ separate the points of $\ga$, and thus 
$\pi_\a:=\oplus_{\f\in\cs_+(\ga)}\pi_\f$
is a faithful representation of $\ga$.
\end{lem}
\begin{proof}
It is enough to see that for $0\neq a\in\ga$, there exists an even state $\om_a$ such that $\om_a(a^*a)>0$. Since the conditional expectation $E=\frac{\idd_\ga+\a}2$ is faithful, there exists $\f\in\cs(\ga_+)$ such that $\f_a(E(a^*a))>0$. It is enough to choose $\om=\f\circ E$.\footnote{If $0\neq a\in\ga$, then either $a_+$ or $a_-$ should be non null, and thus
$E(a^*a)=a_+^*a_++ a_-^*a_->0$.} 
\end{proof}
\medskip

\noindent
\textbf{Fermi algebras.} 
In \cite{CDF}, it has been systematically studied the $C^*$-Fermi systems, that is $C^*$-algebras equipped with a $\bz^2$-grading, and the corresponding $C^*$-Fermi tensor product was constructed by using the maximal norm. In such a context, it was defined and studied the
the maximally entangled state on such a $C^*$-Fermi tensor product $\ga\, \circled{\rm{{\tiny F}}}_{\rm max}\,\ga^\circ$ associated to a single even state 
$\f$ on $\ga$, that is the analogous to the diagonal state associated to $\f$, among those whose marginals are precisely $\f$ on $\ga$, and $\f^\circ$ on the opposite algebra $\ga^\circ$.

In that paper, it was also firstly defined a norm on the $C^*$-Fermi tensor product by using an "universal" representation associated to all even product states, by suggesting that such a norm would  have been the analogous of the minimal $C^*$-cross norm of the usual tensor product.
Later on, in \cite{CRZ} it was shown that such a norm is indeed minimal among all $C^*$-cross norms on the (algebraic) Fermi product of two 
$\bz^2$-graded $C^*$-algebras.

It should be remarked that the minimality of the norm, firstly introduced at the end of Section 8 of \cite{CDF} and used to built the infinite $C^*$-Fermi tensor product, plays no role in the investigation of symmetric states on the infinite $C^*$-Fermi tensor product, but only the fact that it is associated with the faithful representation built by using even product states is relevant for such an analysis.

By coming back to our study, in \cite{Fde} we indeed built the infinite $C^*$-Fermi tensor product $\big(\circled{\rm{{\tiny F}}}_{\,\bn}\, \gb, \circled{\rm{{\tiny F}}}_{\,\bn}\, \a\big)$ of the single $\bz^2$-graded $C^*$-algebra $\gb$ by using such a "minimal" norm. Then we studied the early properties of symmetric states, that is the states in the simplex $\cs_\bp(\ga_{\rm F})$ consisting of all states invariant under the natural action of the group $\bp_\bn=:\bp$ of all finite permutations of 
$\bn$, on $\ga_{\rm F}$. In particular, we extended the quantum De Finetti theorem in \cite{St2} and \cite{CF} to such a general situation. Since only the minimal norm is considered, when we speak about  $C^*$-Fermi tensor products and the usual infinite tensor products, we suppose without further mention, that they are constructed by using the minimal norm.

Let $(\ga,\a)$ and $(\gb,\b)$ be two $\bz^2$-graded $C^*$-algebras.
The key point to construct such a framework is the GNS representation of the product state $\om\times\f$ on 
$\ga\, \circled{\rm{{\tiny F}}}\,\gb$ for even states $\om\in\cs_+(\ga)$ and $\f\in\cs_+(\gb)$. Although the detailed knowledge of such a GNS construction plays no role in the analysis in \cite{Fde}, it is relevant here to prove that the Klein transformation we are building in the foregoing section is realising a $*$-isomorphism between the Fermi tensor product and the usual tensor product.
Such a detailed analysis is carried out in \cite{CRZ}, Section 3, which we report here for the convenience of the reader.\footnote{It should be noticed that, in the more general context in which only one state is even, it was already proved in \cite{CDF}, Section 9, that such a GNS representation of a product state can be constructed by using the Stinespring dilation.}

Indeed, for $\bz^2$-graded $C^*$-algebras $(\ga,\a)$ and $(\gb,\b)$ and even states $\om$ and $\f$ on $\ga$ and $\gb$ as above,
let $(\ch_\om,\pi_\om, V_\om,\xi_\om)$ and $(\ch_\f,\pi_\f, V_\f,\xi_\f)$ be the corresponding covariant (w.r.t. the grading) GNS representations. Since $\bz^2=\{1,-1\}=\{+,-\}$, we are simply denoting by $V_\om:=U_\om(-1)$
and $V_\f:=U_\f(-1)$, the unitaries which implement  $\a$ and $\b$, respectively.
Therefore, $\ch_\om$, $\ch_\f$, and consequently $\cb(\ch_\om)$, $\cb(\ch_\f)$, are graded Banach spaces, see \cite{CDF}, Section 5, and thus it is natural to extend the combined grading to products of operators with vectors. 

For the general situation, fix now a selfadjoint unitary operator $U$ acting on the Hilbert space $\ch$. Since $U=E_+-E_-$, $E_\pm$ being the projections onto the eigenspaces relative to eigenvalues $\pm1$, on one hand $\ch=\ch_-\oplus\ch_-$, and also $\cb(\ch)=\cb(\ch)_+\oplus\cb(\ch)_-$ is endowed with the grading 
induced by $\ad_U$. 
With ''$\partial$'' denoting the grade, it is now meaningful consider the combined grade of operators acting on vectors as follows. 

First we note that, for $T\in\cb(\ch)$ and $\xi\in\ch$ both homogeneous ({\it i.e.} of fixed degree), $T\xi$ is still homogeneous and we have
$\partial(T\xi)=\partial(T)\partial(\xi)$. Second, if we have two graded Hilbert spaces $(\ch,U)$ and $(\ck,V)$, we can  
extend (6.1) in \cite{CDF}, by defining the combined grade for homogeneous elements as
\begin{align*}
\eeps(T,\xi):=\left\{\!\!\!\begin{array}{ll}
                      -1 &\text{if}\,\, \partial(T),\partial(\xi)=-1\,,\\
                     \,\,\,\,\,1 &\text{otherwise}\,,
                    \end{array}
                    \right.
\quad  T\in\cb(\ck)\,,\,\,\xi\in\ch\,.                
\end{align*}

Now we are in position to define the Fermi product $S\, \circled{\rm{{\tiny F}}}\,T$ as follows. Indeed, for $S\in\cb(\ch)$ arbitrary and $T\in\cb(\ck)$ homogeneous, and for $\xi\in\ch$ homogeneous and $\eta\in\ck$ arbitrary, put
\begin{equation}
\label{fermico}
(S\odot T)(\xi\otimes \eta):=\eeps(T,\xi)S\xi\otimes T\eta\,.
\end{equation}

The linear extension of \eqref{fermico} to all elements $S\in\cb(\ch)$, $T\in\cb(\ck)$, and $\xi\in\ch$, $\eta\in\ck$, defines a bounded operator called
$S\, \circled{\rm{{\tiny F}}}\,T\in\cb(\ch\otimes\ck)$.\footnote{With $\ch\otimes\ck$, we denotes the usual Hilbertian tensor product. This differs from the notation $\ch\, \circled{\rm{{\tiny F}}}\,\ck$ in \cite{CRZ}. Since the, perhaps crucial, \eqref{fermico} is not frequently used in the forthcoming section, we prefer to keep the original notation not to weigh down too much the notations. Another evident motivation to keep the original
notation $\ch\otimes\ck$ is that in \eqref{fermico}, the Hilbert space is still the tensor product, and only the combined action of $S$ and $T$ on 
$\ch\otimes\ck$ must change according to the Fermi rule, which of course takes into account also of the grade of elements of $\ch\otimes\ck$ on which they are acting.} 
 
As a consequence, we have the following 
\begin{thm}[\cite{CRZ}, Proposition 3.4]
\label{crz}
The GNS covariant representation $(\ch_{\om\times\f},\pi_{\om\times\f}, V_{\om\times\f},\xi_{\om\times\f})$ of the product state is unitarily equivalent to
$\big(\ch_{\om}\otimes\ck_{\f},\pi_{\om}\, \circled{\rm{{\tiny F}}}\,\pi_{\f},V_\om \otimes V_{\f},\xi_{\om}\otimes\xi_{\f}\big)$.
\end{thm}
\noindent
Notice that, since $V_\om$ and $V_{\f}$ are even, $V_\om\, \circled{\rm{{\tiny F}}}\,V_\f=V_\om\otimes V_\f$.

In addition, it is also possible to define the Fermi-von Neumann tensor product of $bz^2$-graded von Neumann algebras (see \cite{CRZ2}) which takes into account of the $\bz^2$-grading, analogous to the usual von Neumann tensor product, to which it reduces if the grading is trivial. We report here such a definition for the convenience of the reader.

Indeed, for $S\in\cam\subset\cb(\ch)$ and $T\in\cn\subset\cb(\ck)$, $\cam,\cn$ being $\bz^2$-graded von Neumann algebras as above, and $S\,\circled{\rm{{\tiny F}}}\,T$ defined according the rule in \eqref{fermico}, such a Fermi-von Neumann tensor product is simply defined as
$$
\cam\, \overline{\circled{\rm{{\tiny F}}}}\, \cn:=\big\{S\,\circled{\rm{{\tiny F}}}\,T\mid S\in\cam,\,\,T\in\cn\big\}''\,.
$$
 
\medskip

\noindent
\textbf{Infinite $C^*$-Fermi tensor product.} 

We end by recalling that, for a $\bz^2$-graded $C^*$-algebra $(\gb,\b)$ and an even state $\f\in\cs_+(\gb)$, we have a consistent system
$$
\underbrace{\gb\, \circled{\rm{{\tiny F}}}\, \cdots \,\circled{\rm{{\tiny F}}}\,\gb}_{\textrm{(n+1)-times}} \circled{\rm{{\tiny F}}}\,\idd_\gb
\subset
\underbrace{\gb\, \circled{\rm{{\tiny F}}}\, \cdots \, \circled{\rm{{\tiny F}}}\,\gb}_{\textrm{(n+2)-times}}\,,
\quad n\in\bn\,,
$$
of $\bz^2$-graded $C^*$-algebras, compatible with the gradation
$$
\circled{\rm{{\tiny F}}}{\,}_{k=1}^{n+1}\b=\big(\circled{\rm{{\tiny F}}}{\,}_{k=1}^{n+2}\b\big)
\lceil_{\big(\circled{\rm{{\tiny F}}}{\,}_{k=1}^{n+1}\gb\big)\, \circled{\rm{{\tiny F}}}\,\idd_\ga}\,,\quad n\in\bn\,,
$$
and with forming the product states
$$
\times_{k=1}^{n+1}\f=\times_{k=1}^{n+2}\f
\lceil_{\big(\circled{\rm{{\tiny F}}}{\,}_{k=1}^{n+1}\gb\big)\, \circled{\rm{{\tiny F}}}\,\idd_\ga}\,,\quad n\in\bn\,.
$$

Therefore, the $C^*$-inductive limit ({\it e.g.} \cite{WO}) of the direct system as above provides a $\bz^2$-graded $C^*$-algebra
$$
(\ga_{\rm F},\a):=\big(\circled{\rm{{\tiny F}}}_{\,\bn}\, \gb, \circled{\rm{{\tiny F}}}_{\,\bn}\, \b\big)
$$
together with the infinite product state $\times_{\bn}\f\in\cs_+(\ga_{\rm F})$ of the single state $\f\in\cs_+(\gb)$. The reader is referred to \cite{Fde} for further details.

\section{The Klein transformation}
\label{flhf}

The present section is devoted to the extension of the Klein transformation reported in \cite{T}, Exercise XIV.1, to quite  general situations. We limit the analysis to the situation relative to the innerness of the grading of one of the involved algebras, the remaining situations being not coverable as explained in Section \ref{ezx}.

Indeed, we start with
$\bz_2$-graded $C^*$-algebra $(\ga,\a)$, $(\gb,\b)$, and suppose that the grading of one of them, say that of $\ga$, is generated by an inner unitary, that is
$\a=\ad_u$, for some selfadjoint (necessarily even) unitary $u\in\ga$. This would cover many cases of interest, including when $\ga$ is a type $\ty{I}$ factor, and in particular the cases of full matrix algebras, including the CAR, the last being the only case managed in some detail up to now ({\it e.g.} \cite{AFM, AM1, CF, CF1, CDF, F1}).
\begin{prop}
\label{klein}
Let $(\ga,\a)$, $(\gb,\b)$ $\bz_2$-graded $C^*$-algebras, with $\a=\ad_u$ for some selfadjoint unitary $u\in\ga$. Then the $\bz_2$-Fermi $C ^*$-tensor product $(\ga\, \circled{\rm{{\tiny F}}}\, \gb,\a\, \circled{\rm{{\tiny F}}}\,\b)$ is $*$-isomorphic to the $\bz_2$-tensor product $(\ga\, \circled{\rm{{\tiny X}}}\, \gb,\a\, \circled{\rm{{\tiny X}}}\,\b)$.
\end{prop}
\begin{proof}
We start by defining
 \begin{equation}
 \label{klevo0}
\ga\,\circled{\rm{{\tiny F}}}\,\gb \ni a\,\circled{\rm{{\tiny F}}}\,b
\mapsto\s(a\,\circled{\rm{{\tiny F}}}\,b)
:=a\,\circled{\rm{{\tiny F}}}\,b_++(au)\,\circled{\rm{{\tiny F}}}\,b_-\in\ga\,\circled{\rm{{\tiny F}}}\,\gb\,.
\end{equation}

We first note that $\s$ is a real map. Indeed, for the 2-nd piece in the r.h.s. above,
\begin{align*}
\big((au)\,\circled{\rm{{\tiny F}}}\,b_-\big)^*=&\big(\idd_\ga\,\circled{\rm{{\tiny F}}}\,b_-^*\big)\big((au)^*\,\circled{\rm{{\tiny F}}}\,\idd_\gb\big)\\
=&\big(\idd_\ga\,\circled{\rm{{\tiny F}}}\,b_-^*\big)\big((ua^*_+-ua^*_-)\,\circled{\rm{{\tiny F}}}\,\idd_\gb\big)\\
=&(ua^*_++ua^*_-)\,\circled{\rm{{\tiny F}}}\,b_-=\big((a^*u)\,\circled{\rm{{\tiny F}}}\,b^*_-\big)\,.
\end{align*}

It is also grading equivariant, which is immediately achieved because the unitary $u$ is even. Since the checking relative to equivariance relative to the involved gradings is elementary, we omit to discuss further this aspect. Instead, it is not a homomorphism, provided it is 
well-defined.\footnote{That $\s$ is not a homomorphism is natural, and emerges already in the simplest case of the CAR algebra, see {\it e.g.} \cite{T}, Exercise
XIV.1.}

We check the claimed statements, and in particular that $\s$ is indeed well-defined. For such a purpose we consider the universal representation
$$
\pi=\stackrel[\f\in\cs_+(\ga),
\psi\in\cs_+(\gb)]{}{\bigoplus}\pi_{\f\times\psi}
=\stackrel[\f\in\cs_+(\ga),\psi\in\cs_+(\gb)]{}{\bigoplus}\pi_{\f}\, \circled{\rm{{\tiny F}}}\, \pi_{\psi}\,,
$$
acting on the direct sum of  $\ch_\f\otimes\ch_\psi$, whereas the gradings are implemented by the corresponding direct sum of
$\pi_\f(u)\otimes V_\psi$.\footnote{Since $\pi_\f(u)$ and $V_\psi$ are both even, their Fermi tensor product $\pi_\f(u)\,\circled{\rm{{\tiny F}}}\,V_\psi$, defined according the rule in \cite{CRZ} and reported in Section \ref{sec1}, coincides with the usual tensor product  
$\pi_\f(u)\,\circled{\rm{{\tiny X}}}\,V_\psi\equiv\pi_\f(u)\otimes V_\psi$.}

For $a\in\ga$ and $b\in\gb$, put 
\begin{equation}
\label{abi}
\begin{split}
&A:=a\, \circled{\rm{{\tiny F}}}\, \idd_{\gb}\,,\\
&B:=\s(\idd_{\ga}\,\circled{\rm{{\tiny F}}}\, b)=
\idd_{\ga}\, \circled{\rm{{\tiny F}}}\,b_+
+u\, \circled{\rm{{\tiny F}}}\,b_-\,,
\end{split}
\end{equation}
and compute
\begin{align*}
AB=&a\, \circled{\rm{{\tiny F}}}\,b_+
+(au)\, \circled{\rm{{\tiny F}}}\,b_-\,,\\
BA=&a\, \circled{\rm{{\tiny F}}}\,b_++(ua_+)\, \circled{\rm{{\tiny F}}}\,b_-
-(ua_-)\, \circled{\rm{{\tiny F}}}\,b_-\\
=&a\, \circled{\rm{{\tiny F}}}\,b_++(a_+u)\, \circled{\rm{{\tiny F}}}\,b_-
+(a_-u)\, \circled{\rm{{\tiny F}}}\,b_-\\
=&a\, \circled{\rm{{\tiny F}}}\,b_+
+(au)\, \circled{\rm{{\tiny F}}}\,b_-=AB\,.
\end{align*}
 
Therefore, the generators $A$, $B$ commute each other, which leads to 
$\big[\ga\, \circled{\rm{{\tiny F}}}\, \idd_{\gb},\s(\idd_{\ga}\,\circled{\rm{{\tiny F}}}\, \gb)\big]=0$.

Now, for $A,B$ as above, a product state $\om\times\f$, $\om\in\cs_+(\ga)$, $\f\in\cs_+(\gb)$, 
$\xi,x\in\ch_\om$, $\eta,y\in\ch_\f$ and, by taking into account \eqref{fermico}, we compute
\begin{align*}
\pi_{\om\times\f}(AB)(\xi\otimes\eta)=&\big(\pi_\om(a)\, \circled{\rm{{\tiny F}}}\, \pi_\f(b_+)+(\pi_\om(a)V_\om)\, \circled{\rm{{\tiny F}}}\, \pi_\f(b_-)\big)\xi\otimes\eta\\
=&\pi_\om(a)\xi\otimes\pi_\f(b_+)\eta+\big(\pi_\om(a)V_\om\xi_+\otimes\pi_\f(b_+)\eta\\
-&\pi_\om(a)V_\om\xi_+\otimes\pi_\f(b_+)\eta\big)=\pi_\om(a)\xi\otimes\pi_\f(b_+)\eta\\
+&\big(\pi_\om(a)\xi_+\otimes\pi_\f(b_+)\eta
+\pi_\om(a)\xi_+\otimes\pi_\f(b_+)\eta\big)\\
=&\pi_\om(a)\xi\otimes\pi_\f(b)\eta\,.
\end{align*}

Finally, taking first finite linear combinations of generators $A$ and $B$, then of vectors in $\ch_{\om\times\f}=\ch_{\om}\otimes\ch_{\f}$ such that
$$
\Big\|\sum_{\xi,\eta}\xi\otimes \eta\Big\|=1=\Big\|\sum_{x,y}x\otimes y \Big\|\,,                
$$
and finally taking the supremum for any fixed combination of generators, we get after the last computation,
\begin{align*}
\Big\|\sum_{A,B}AB\Big\|=&\sup_{\om,\f}\sup_{\xi,\eta;x,y}
\Big|\sum_{\xi,\eta,x,y}\sum_{a,b}\big\langle\pi_\om(a)\xi\otimes\pi_\f(b)\eta,x\otimes y\big\rangle\Big|\\
=&\sup_{\om,\f}\sup_{\xi,\eta;x,y}
\Big|\sum_{\xi,\eta,x,y}\sum_{a,b}\big\langle\pi_\om(a)\xi\otimes\pi_\f(b)\eta,x\otimes y\big\rangle\Big|\\
=&\sup_{\om,\f}\sup_{\xi,\eta;x,y}\Big|\sum_{a,b}\Big(\sum_{\xi,x}\langle\pi_\om(a)\xi,x\rangle\Big) 
\Big(\sum_{\eta,y}\langle\pi_\f(b)\eta,y\rangle\Big) \Big|\\
=&\Big\|\sum_{a,b}a\, \circled{\rm{{\tiny X}}}\,b\Big\|\,.
\end{align*}

Since, $\pi$, $\pi_\a$ and $\pi_\b$ are faithful, $\s$ given in \eqref{klevo0} is a well-defined map of $\ga\odot\gb$ which extends to an 
invertible bicontinuous map of $\ga\,\circled{\rm{{\tiny F}}}\,\gb$ onto itself, and in addition $\ga\,\circled{\rm{{\tiny F}}}\,\gb$ is generated by 
$\ga\, \circled{\rm{{\tiny F}}}\, \idd_{\gb}\sim\ga$ and $\s(\idd_{\ga}\,\circled{\rm{{\tiny F}}}\, \gb)\sim\gb$. Moreover, such copies of $\ga$ and $\gb$ inside $\ga\,\circled{\rm{{\tiny F}}}\,\gb$ generate a tensor product :
$$
*-{\rm alg}\big\{(\ga)\, \circled{\rm{{\tiny F}}}\, \idd_{\gb}),
\s\big(\idd_{\ga}\, \circled{\rm{{\tiny F}}}\,\gb)\big\}\sim\ga\otimes\gb\,,
$$
where ''$\sim$'' means isometric $*$-isomorphism.

Summarising, we get
\begin{align*}
(\ga\, \circled{\rm{{\tiny F}}}\, \gb,\a\, \circled{\rm{{\tiny F}}}\,\b)
=&\Big(\overline{*-{\rm alg}\big\{(\ga)\, \circled{\rm{{\tiny F}}}\, \idd_{\gb}),
\s\big(\idd_{\ga}\, \circled{\rm{{\tiny F}}}\,\gb)\big\}},
\a\, \circled{\rm{{\tiny F}}}\,\b\Big)\\
\sim&(\ga\, \circled{\rm{{\tiny X}}}\, \gb,\a\, \circled{\rm{{\tiny X}}}\,\b)\,.
\end{align*}
\end{proof}

Now we go into details concerning such a fundamental identification, and show that, in addition to preserving the grades, it preserves also the product states.
For such a purpose we compute $\om\times\f$ on the generators given in \eqref{abi}, obtaining
\begin{align*}
\om\times\f(AB)=&\om\times\f(a_+\, \circled{\rm{{\tiny F}}}\,b_+)+\om\times\f(a_-\, \circled{\rm{{\tiny F}}}\,b_+)\\
+&\om\times\f((a_+u)\, \circled{\rm{{\tiny F}}}\,b_-)+\om\times\f((a_-u)\, \circled{\rm{{\tiny F}}}\,b_-)\\
=&\om(a_+)\f(b_+)=\om(a)\f(b)\,.
\end{align*}

Let now $\k_{\ga,\gb}: \ga\, \circled{\rm{{\tiny F}}}\, \gb\to \ga\, \circled{\rm{{\tiny X}}}\, \gb$ realising the equivalence between 
$\ga\, \circled{\rm{{\tiny F}}}\, \gb$ and $\ga\, \circled{\rm{{\tiny X}}}\, \gb$. It is grading equivariant and,
by the above calculations, we see that $\k^{\rm t}_{\ga,\gb}(\psi_{\om,\f})=\om\times\f$ for every states $\om\in\cs(\ga)_+$ and $\f\in\cs(\gb)_+$. 
Such an isomorphism is called a {\it Klein transformation}.

We then have proven the following
\begin{thm}
\label{klcaz}
The Klein transformation $\k_{\ga,\gb}:\ga\, \circled{\rm{{\tiny F}}}\, \gb\to \ga\, \circled{\rm{{\tiny X}}}\, \gb$ 
realises the $*$-isomorphism between $\ga\, \circled{\rm{{\tiny F}}}\, \gb$ and $\ga\, \circled{\rm{{\tiny X}}}\, \gb$, whose transpose 
$\k_{\ga,\gb}^{\rm t}$ provides a bijective map between 
even product states on $\ga\, \circled{\rm{{\tiny X}}}\, \gb$ and (even) product states on $\ga\, \circled{\rm{{\tiny F}}}\, \gb$.
\end{thm}
\begin{proof}
It follows collecting together Proposition \ref{klein} and the previous considerations.
\end{proof}
\begin{rem}
Notice that the Klein transformation does not preserve the slices. However, 
$\k_{\ga,\gb}(a_+\, \circled{\rm{{\tiny F}}}\, b_+)=(a_+\, \circled{\rm{{\tiny X}}}\, b_+)$, which is enough to exhaustively manage the even product states.
\end{rem}

In view to the applications to the investigations of some relevant properties of the symmetric states, we extend the Klein transformation to the infinite tensor products.

Indeed, first we note that the grading of 
$\underbrace{\gb\, \circled{\rm{{\tiny F}}}\, \cdots \, \circled{\rm{{\tiny F}}}\,\gb}_{\textrm{(n+1)-times}}\sim
\underbrace{\gb\, \circled{\rm{{\tiny X}}}\, \cdots \, \circled{\rm{{\tiny X}}}\,\gb}_{\textrm{(n+1)-times}}$, is innerly generated for each $n=1,2,\dots$, and thus one can perform the inductive process to construct  $\k_n:=\k_{\,\circled{\rm{{\tiny F}}}_{k=0}^n\,\ga,\ga}$ for $n=1,2,\dots$\,.\footnote{The reader is referred to \cite{T}, Exercise XIV.1, to recognise how such an inductive process is working for the simplest case of ${\rm CAR}(\bn)$.}

It is now matter of routine to check that such a sequence of Klein transformations $( k_n)_n$
\begin{itemize}
\item[(i)] preserves the corresponding gradings: 
$$
\k_n\circ\big(\,\circled{\rm{{\tiny F}}}_{k=0}^n\a\big)=\big(\,\circled{\rm{{\tiny X}}}_{k=0}^n\a\big)\circ\k_n\,;
$$
\item[(ii)] sends product states made of a single even one in the corresponding product state: 
$$
\k^{\rm t}_n\big(\psi_{\underbrace{\f,\f,\dots}_{\textrm{(n+1)-times}}}\big)
=\underbrace{\f\times\f\times\cdots}_{\textrm{(n+1)-times}}\,,\quad \f\in\cs_+(\gb)\,;
$$
\item[(iii)] is compatible:
$$
\k_{n+1}(a\,\circled{\rm{{\tiny F}}}\, \idd)=
 \k_n(a)\,\circled{\rm{{\tiny F}}}\, \idd, \,\,n=1,2,\dots,\,\,\,
 a\in\underbrace{\gb\, \circled{\rm{{\tiny F}}}\, \cdots \, \circled{\rm{{\tiny F}}}\,\gb}_{\textrm{(n+1)-times}}\,.
$$
\end{itemize}
 
Therefore, we have the following
\begin{thm}
\label{klete}
Let $(\gb,\b)$ be a $\bz_2$-graded $C^*$-algebra. The $C^*$-inductive limit $\,\lim_{n\to\infty}\k_n=:\k$
of the sequence $\k_n$, $n=1,2,\dots$, exists 
and 
\begin{itemize}
\item[(a)] realises a $*$-isomorphism between $\ga_{\rm F}$ and $\ga_{\rm X}$; 
\item[(b)] preserves the corresponding gradings: 
$$
\k\circ\big(\,\circled{\rm{{\tiny F}}}_{\bn}\a\big)=\big(\,\circled{\rm{{\tiny X}}}_{\bn}\a\big)\circ\k\,;
$$
\item[(c)] sends product states made of a single even state in $\cs(\ga_{\rm X})$ in the corresponding one in $\cs(\ga_{\rm F})$:
$$
\k^{\rm t}\big(\psi_{\f,\f,\dots}\big)
=\times_\bn\f\,,\quad\f\in\cs_+(\gb)\,.
$$
\end{itemize}
\end{thm}
\begin{proof}
The assertion follows by \cite{WO}, Appendix L.1, by taking into account (i)-(iii) above.
\end{proof}

\section{The Klein transformation and the symmetric states}

Here we prove a result which, on one hand can be applied to the investigation of some relevant properties of symmetric states, on the other hand has a self containing interest. For such a purpose, we denote a Choquet simplex simply by a simplex.

Such an analysis can carried out only when the grading automorphism $\b=\ad_u$ of the $\bz_2$-graded $C^*$-algebra $\gb$ is generated by an inner selfadjoint uniatry $u\in\gb$, which we assume throughout the present section without any further mention.
\begin{thm}
\label{symis}
Let $(\gb,\b)$ be $\bz_2$-graded $C^*$-algebra. Then the restriction of transposed $\k^{\rm t}$ of Klein transformation in Theorem \ref{klete}, induces topological isomorphisms between simplexes
\begin{align*}
\cs_\bp\big(\,\circled{\rm{{\tiny X}}}_{\bn}\gb_+\big)=\cs_\bp\big(\,\circled{\rm{{\tiny F}}}_{\bn}\gb_+\big)
\sim&\cs_\bp\big(\,\circled{\rm{{\tiny F}}}_{\bn}\gb\big)\\
=\cs_{\bp\times\bz^2}\big(\,\circled{\rm{{\tiny F}}}_{\bn}\gb\big)
=&\k^{\rm t}\Big(\cs_{\bp\times\bz^2}\big(\,\circled{\rm{{\tiny X}}}_{\bn}\gb\big)\Big)\,.
\end{align*}
\end{thm}
\begin{proof}
The first congruence can be proved by using the transposed of the canonical conditional expectation 
$E:\,\circled{\rm{{\tiny F}}}_{\bn}\gb\to\,\circled{\rm{{\tiny F}}}_{\bn}\gb_+=\circled{\rm{{\tiny X}}}_{\bn}\gb_+$ and is left to the 
reader.\footnote{See \cite{CDF}, Section 9, for some detail about this point.} We prove only that induced by the Klein transformation.

First we note that the $C^*$-dynamical system $\big(\,\circled{\rm{{\tiny X}}}_{\bn}\gb,\bp\times\bz^2\big)$ is asymptotically abelian in norm, for it is enough to do the same computations in \cite{St2} by considering the sequence $(g_n\times e_{\bz^2})_n\subset\bp\times\bz^2$, where 
$(g_n)_n\subset\bp$ is given in \cite{St2}, Formula (1). Therefore $\big(\,\circled{\rm{{\tiny X}}}_{\bn}\gb,\bp\times\bz^2\big)$ is
$\bp\times\bz^2$-abelian, and thus $\cs_{\bp\times\bz^2}\big(\,\circled{\rm{{\tiny X}}}_{\bn}\gb\big)$ is a simplex.

Now $\k^{\rm t}$, which is a bijective and bicontinuous, maps $\ce_{\bp\times\bz^2}\big(\,\circled{\rm{{\tiny X}}}_{\bn}\gb\big)$ onto
$\ce_{\bp\times\bz^2}\big(\,\circled{\rm{{\tiny F}}}_{\bn}\gb\big)=\ce_{\bp}\big(\,\circled{\rm{{\tiny F}}}_{\bn}\gb\big)$, and thus it maps
$\cs_{\bp\times\bz^2}\big(\,\circled{\rm{{\tiny X}}}_{\bn}\gb\big)$ onto $\cs_{\bp}\big(\,\circled{\rm{{\tiny F}}}_{\bn}\gb\big)$ by the Krein-Milman theorem.
Namely, $\k^{\rm t}\lceil_{\cs_{\bp\times\bz^2}\big(\,\circled{\rm{{\tiny X}}}_{\bn}\gb\big)}$ is an affine map which realises a topological homeomorphism between the closed convex sets 
$\cs_{\bp\times\bz^2}\big(\,\circled{\rm{{\tiny X}}}_{\bn}\gb\big)$ and $\cs_{\bp}\big(\,\circled{\rm{{\tiny F}}}_{\bn}\gb\big)$.\footnote{Here, affine means that for the segment $\l\in[0,1]$ and $\om,\f\in\cs_{\bp\times\bz^2}\big(\,\circled{\rm{{\tiny X}}}_{\bn}\gb\big)$, 
$\k^{\rm t}\big(\l\om+(1-\l)\f\big)=\l\k^{\rm t}(\om)+(1-\l)\k^{\rm t}(\f)$).}

It remains to see that it preserves the fact to being simplexes,  
that is the restriction of $\k^{\rm t}$ to positive $\bp\times\bz^2$-invariant functionals on $\,\circled{\rm{{\tiny X}}}_{\bn}\gb$
preserves its lattice structure. For we proceed as in the proof of Theorem 3.1.14 of \cite{S}.\footnote{The fact that restriction of $\k^{\rm t}$ realises a topological isomorphisms of simplexes would follow by \cite{BR1}, Theorem 4.1.15. For the reader convenience, we are providing a direct proof.}

First we note that the GNS covariant representation of $\f\in\cs_{\bp}(\ga_{\rm F})$, with $\f=\psi\circ\k$ and $\psi\in\cs_{\bp\times\bz^2}(\ga_{\rm X})$, is 
$\big(\ch_\psi,\pi_\psi\circ\k, V_\psi,\xi_\psi\big)$. 

Fix now two positive (non null) functionals $f_1,f_2\in\ga_{\rm X}^*$ invariant under the action of $\bp\times\bz^2$, and define $\psi:=\frac{f_1+f_2}{\|f_1\|+\|f_2\|}$. Then $\psi\in\cs_{\bp\times\bz^2}(\ga_{\rm X})$, and
$$
f_i(a)=\langle\pi_\psi(a)h_i\xi_\psi,\xi_\psi\rangle\,,\quad a\in\ga_{\rm X},\,\,i=1,2\,,
$$
where the $h_i$ are uniquely determined densities in the selfadjoint part of the abelian algebra $\{\pi_\psi(\ga_{\rm X}),V_\psi(\bp\times\bz^2)\}'$. Therefore,
$$
f_1\vee f_2(a)=\langle\pi_\psi(a)(h_1\vee h_2)\xi_\psi,\xi_\psi\rangle\,,\quad a\in\ga_{\rm X},\,\,i=1,2\,,
$$
and thus
$$
\big(\k^{\rm t}(f_1\vee f_2)\big)(a)=(f_1\vee f_2)(\k(a))
=\langle\pi_\psi(\k(a))(h_1\vee h_2)\xi_\psi,\xi_\psi\rangle\,.
$$

On the other hand, reasoning as above, we have 
$$
\k^{\rm t}f_i(a)=\langle\pi_\psi(\k(a))h_i\xi_\psi,\xi_\psi\rangle\,,\quad a\in\ga_{\rm F},\,\,i=1,2\,,
$$
and thus 
$$
(\k^{\rm t}f_1\vee \k^{\rm t}f_2)(a)=\langle\pi_\psi(\k(a))(h_1\vee h_2)\xi_\psi,\xi_\psi\rangle\,,\quad a\in\ga_{\rm F}\,,
$$

Collecting together, we have $\k^{\rm t}(f_1\vee f_2)=\k^{\rm t}f_1\vee \k^{\rm t}f_2$ which is the assertion.
\end{proof}

As an immediate consequence, we immediately provides the proof of the analogous result in \cite{St2}, Theorem 2.8, and reported also in \cite{Fde} as Proposition 6.3. However, such a result holds true without any further assumption, in particular without assume the innerness of the grading automorphism).

Then we sketch the proof of the following
\begin{thm}
\label{boundary} 
The Choquet simplex $\cs_{\bp}(\ga_{\rm F})$
has a $*$-weakly closed boundary and is affinely isomorphic to
the probability Radom measures on the convex $*$-weakly compact set $\cs(\gb_+)$.
\end{thm}
\begin{proof}
(sketch) Let $E: \gb\to\gb_+$. It is not hard to convince oneself that, analogously to 
$$
E_{\rm X}:=\circled{\rm{{\tiny X}}}_{\bn}\,E:\,\circled{\rm{{\tiny X}}}_{\bn}\gb\to\,\circled{\rm{{\tiny X}}}_{\bn}\gb_+\,,
$$
one can construct 
$$
E_{\rm F}:=\circled{\rm{{\tiny F}}}_{\bn}\,E:\,\circled{\rm{{\tiny F}}}_{\bn}\gb\to\,\circled{\rm{{\tiny X}}}_{\bn}\gb_+\,.
$$

It is also not hard to see that $E_{\rm X}^{\rm t}$ and $E_{\rm F}^{\rm t}$ realise isomorphisms of simplexes 
$\cs_\bp\big(\,\circled{\rm{{\tiny X}}}_{\bn}\gb_+\big)\sim\cs_{\bp\times\bz^2}\big(\,\circled{\rm{{\tiny X}}}_{\bn}\gb\big)$ and
$\cs_\bp\big(\,\circled{\rm{{\tiny F}}}_{\bn}\gb_+\big)\sim\cs_\bp\big(\,\circled{\rm{{\tiny F}}}_{\bn}\gb\big)$, respectively.

Now, since $\circled{\rm{{\tiny F}}}_{\bn}\gb_+=\,\circled{\rm{{\tiny X}}}_{\bn}\gb_+$, 
we conclude that 
$\cs_{\bp\times\bz^2}\big(\,\circled{\rm{{\tiny X}}}_{\bn}\gb\big)\sim
\cs_{\bp}\big(\,\circled{\rm{{\tiny F}}}_{\bn}\gb\big)$. 

The claim follows because $\cs_{\bp\times\bz^2}(\ga_{\rm X})$ is affinely isomorphic with the probability Radom measures on the convex $*$-weakly compact set 
$\cs(\gb_+)$.
\end{proof}

\section{The simplex of symmetric states}

In this section, we study the structure and the properties of the set of symmetric states on Fermi systems. As we have seen in Section \ref{ezx}, it is not possible to directly compare the set of the symmetric states on $\ga_{\rm F}$ with the analogous one on $\ga_{\rm X}$, even if they share many common properties.

In the first part, we prove such shared properties 
whereas in the second part, assuming that the grading is innerly generated, we take advantage of the associated Klein transformation to prove some other relevant properties.

We then start with the results which are directly connected with the fact that $\bp$ acts on $\ga_{\rm F}$, again as a large group of automorphisms (see \cite{St1}), and thus without assuming innerness.

For $\#\in\{\ty{I}, \ty{II_1},\ty{II_\infty}, \ty{III}\}$, let
$$
\cs_\bp(\ga_{\rm F})_\#:=\{\om\in\cs_\bp(\ga)\mid \pi_\om(\ga_{\rm F})''\,\,\text{is of type}\,\,\#\}\,.
$$

Recall that a {\it face} of the simplex $\cs_\bp(\ga_{\rm F})$ is a convex subset $F$ such that, if $\om\in F$ dominates $\f\in\cs_\bp(\ga_{\rm F})$ ({\it i.e.} $\f\leq\l\om$ for some $\l>0$), then also $\f\in F$.

\begin{prop}
\label{face}
For each $\#$, $\cs_\bp(\ga_{\rm F})_\#$ is a face of $\cs_\bp(\ga)$.
\end{prop}
\begin{proof}
It is the same as in Lemma 2.9 of \cite{St2}. We sketch the proof for the convenience of the reader.

Suppose that $\cs_\bp(\ga)\ni\om=\l\om_1+(1-\l)\om_2$, $\l\in(0,1)$, is a non trivial convex combination of $\om_1,\om_2\in\cs_\bp(\ga_{\rm F})_\#$. Since $\bp$ is acting on $\ga$ as a large group of automorphisms ({\it cf.} \cite{Fde}, Proposition 5.3), there exist $B_i\in\ty{Z}_\om^\bp\subset\ty{Z}_\om$ such that 
$\om_i=\om_{A_i\xi_\om}$, $i=1,2$. 

Let $P_i$ be the selfadjoint projection onto the essential subspace $[\pi_\om(\ga)A_i\xi_\om]$. One recognise that, first the $P_i$ are central projections, and then $P:=P_1\bigvee P_2=P_1+(P_2-P_1P_2)=I$,
$\pi_{\om}\cong P_i\pi_{\om_i}$, $i=1,2$.

Therefore, the $P_i\pi_\om(\ga)''$ are of type $\#$ by assumption, and, after cutting the common part associated to $P_1P_2$, also 
$(P_2-P_1P_2)\pi_\om(\ga)''$ is of type $\#$.

Collecting together, since $P_1+P_2=I$ and thus 
$$
\pi_\om(\ga)''=P_1\pi_\om(\ga)''+(P_2-P_1P_2)\pi_\om(\ga)''\,,
$$
we argue $\om\in\cs_\bp(\ga_{\rm F})_\#$.

With the same arguments as before, we can argue that, if $\cs_\bp(\ga)\ni\om\leq \l\f$ for some $\l>0$ with $\f\in\cs_\bp(\ga_{\rm F})_\#$, we find a nonzero projection 
$P\in\ty{Z}_\f^\bp\subset\ty{Z}_\f$ such that $\pi_{\om}\cong P_i\pi_{\f}$, and thus $P\pi_\f(\ga)''$ is of the same type of $\pi_\f(\ga)''$.
\end{proof}
If $\om\in\cs_\bp(\ga_{\rm F})$, denote by $E_\#^\om\in{\rm Z}_\f$ the central projection associated to the type $\#$ portion of $\pi_\om(\ga)''$.
\begin{rem}
Since $E_\#^\om$ is $\ad_{U_\om(g)}$-invariant, $g\in\bp$, we argue as in Corollary 2.10 of \cite{St2} that $\cs_\bp(\ga_{\rm F})$ is the convex hull of its faces
$\cs_\bp(\ga_{\rm F})_\#$,  $\#=\ty{I}, \ty{II_1},\ty{II_\infty}, \ty{III}$.
\end{rem}

Let $\l\in[0,1]$, and $\gam$ be a $W^*$-algebra with separable predual, and consider the central decomposition
$\gam=\int^\oplus_\G \gam_\g\di\n(\g)$, see {\it e.g. }\cite{S}. Here, $\G$ is a locally compact dense subset of $\s({\rm Z}(\gam))$ and $\n$ is a Radon probability measure. By Theorem 21.2 of \cite{N}, if $\l\in[0,1]$, the set 
$$
E_\l:=\{\g\in\G\mid \gam_\g\,\text{is a factor of type}\,\, \ty{III}{\,\!}_\l\}
$$
is measurable, and thus it is meaningful to speak about $W^*$-algebras of type $\ty{III}_\l$ as those for which $\n(E_\l)=1$. 

Notice that, if $p$ is a central projection, $p=\int^\oplus_\G\chi_{A_p}(\g)\di\n(\g)$ where $\chi_{A_p}$ is the indicator-function of the measurable set $A$, which is uniquely determined modulo negligible sets by $p$, and vice-versa. Therefore, 
$$
p\,\gam=\int^\oplus_{A_p}\gam_\g\di\n(\g)\,, \quad p\,\, \text{projection of}\,\, \gam\,,
$$
and thus $\gam$ is of type $ \ty{III}_\l$ if and only if $p\,\gam$ is of type $\ty{III}_\l$ for each nonzero projection $p\in{\rm Z}(\gam)$. 

Let now the $\bz^2$-graded $C^*$-algebra $\gb$ be separable. For each $\l\in[0,1]$, define
$$
\cs_\bp(\ga_{\rm F})_\l:=\{\om\in\cs_\bp(\ga)\mid \pi_\om(\ga_{\rm F})''\,\,\text{is of type}\,\, {\rm III}_\l\} 
$$
\begin{prop}
Suppose that $\gb$ is separable. Then $\cs_\bp(\ga_{\rm F})_\l$ is a face of $\cs_\bp(\ga_{\rm F})$.
\end{prop}
\begin{proof}
We start by noticing that, since $\gb$ and thus $\ga_{\rm F}$, is separable, all $\pi_\om$, $\om\in\cs_\bp(\ga_{\rm F})$ act on a separable Hilbert space. Therefore, it is meaningful to speak about the type for type $\ty{III}$ representations.

Now, we can argue as in Proposition \ref{face}. First we fix a convex combination $\om:=\l\om_1+(1-\l)\om_2$, $\l\in(0,1)$, is convex combination of $\om_1,\om_2\in\cs_\bp(\ga_{\rm F})_\l$. Then we find two central projections $P_1,P_2$, summing up to 1, such that 
$\pi_\om(\ga)''=P_1\pi_\om(\ga)''+(P_2-P_1P_2)\pi_\om(\ga)''$ where the addenda are both of type $\ty{III}_\l$ by the above considerations. 

When we have $\cs_\bp(\ga)\ni\om\leq\m\f$, for $\m>0$ and $\f\in\cs_\bp(\ga_{\rm F})_\l$, we argue at the same time that there exists a central projection in 
$\pi_\f(\ga)''$
such that, after identifying $\pi_\om$ with $P\pi_\f\lceil_{\ch_\om}$, $\pi_\om(\ga)''=P\pi_\om(\ga)''$. Hence, $\om\in\cs_\bp(\ga_{\rm F})_\l$ as well.
\end{proof}
\begin{rem}
We have incidentally proved that, even for the usual tensor product, $\cs_\bp(\ga_{\rm X})_\l$ is a face of $\cs_\bp(\ga_{\rm X})$.
\end{rem}
To conclude this first part, by following the approach in \cite{CF1},
a quantum stochastic process with index-set $\bn$ concerning Fermi systems, can be viewed simply as a state on $\ga_{\rm F}$. The {\it exchangeable processes} correspond to symmetric states (w.r.t the natural action of the permutation group).

In addition, for the extension of the notion of independence to quantum situation, the concept of conditional independence w.r.t. the tail algebra (or the algebra at infinity in physical language) is crucial as in classical probability. The readable definition for a stochastic process to being {\it conditionally independent and identically distributed w.r.t the tail algebra} is provided in \cite{Fbo}, 
Definition 1.
\begin{rem}
By using exactly the same proof as that of Theorem 5.4 in \cite{CF1} relative to the CAR algebra, it would follow that a state $\om\in\cs(\ga_{\rm F})$ is symmetric if and only if the corresponding process is conditionally independent and identically distributed w.r.t the tail algebra.
\end{rem}

Now we pass to the analogous of those results in \cite{St2} whose direct proof requires the use of the Klein transformation. Therefore, for the rest of the present section we assume that, for the $\bz^2$-graded $C^*$-algebra $(\gb,\b)$, $\b=\ad_u$ for some, necessarily even, unitary $u\in\gb$. 

However, we would like to point out that the condition of innerness is automatically satisfied if $\gb$ is a 
type $\ty{I}$ von Neumann factor, covering many interesting cases, including the CAR for which $\gb\sim\bm_2(\bc)$.

\begin{prop}
\label{stscu}
Let $\gb$ be a $C^*$-algebra, and $\f\in\cs_+(\gb)$ an even state. Then the infinite product state $\psi_{\f,\f,\dots}$
on $\ga:=\,\circled{\rm{{\tiny X}}}_{\,\bn}\, \gb$ has central support in the bidual if and only if $\f$ has. The same is true for the corresponding product state $\times_\bn\f$ on $\ga:=\,\circled{\rm{{\tiny F}}}_{\,\bn}\, \gb$.
\end{prop}
\begin{proof}
The second part directly follows by taking into account that, if $\f=\om\circ\k$, then $\pi_\f\cong\pi_\om\circ\k$, and thus 
$\pi_\f(\ga_{\rm F})''\cong\pi_\om(\ga_{\rm X})''$. The first part is quite well-known, but we report that for convenience of the reader.

Indeed, the GNS representation of the infinite product state $\psi_{\f,\f,\dots}$ is acting on the incomplete infinite tensor product $\otimes_\bn\ch_\f$, w.r.t. the infinite constant sequence $\xi_\f,\xi_\f,\dots$, which provides also the cyclic vector:
$\xi_{\psi_{\f,\f,\dots}}=\otimes_\bn\xi_\f$.
The representation is given for the total set of elements of the form $a:=a_0\otimes a_1\otimes\cdots\otimes a_n\otimes\idd_\ga\otimes\idd_\ga\otimes\cdots$, in the usual way by
$$
\pi_{\psi_{\f,\f,\dots}}(a)=\pi_\f(a_0)\otimes\pi_\f(a_1)\otimes\cdots\otimes\pi_\f(a_n)\otimes\idd_{\ch_\f}\otimes\idd_{\ch_\f}\otimes\cdots\,,
$$
see e.g. \cite{vN}, or \cite{St}, A 17.

The commutant $\pi_{\om}(\ga)'$ is generated by the $*$-algebra 
$\pi_{\om}(\ga)'_o$ given by
$$
*-{\rm alg}\bigg\{\underbrace{\pi_\f(\gb)'\ots\pi_\f(\gb)'\ots\cdots\ots\pi_\f(\gb)'}_{\textrm{(n+1)-times}}\otimes\idd_{\ch_\f}\otimes\idd_{\ch_\f}\otimes\cdots
\mid n\in\bn\bigg\}\,.
$$

Suppose that the product state $\psi_{\f,\f,\dots}$ has central support in the bidual. Then $\otimes_\bn\xi_\f$ is cyclic for the above 
weakly dense $*$-algebra, and thus 
\begin{align*}
\big[(\pi_\f(\gb)'\otimes\idd_{\ch_\f}&\otimes\idd_{\ch_\f}\otimes\cdots)(\xi_\f\otimes\xi_\f\otimes\xi_\f\otimes\cdots)\big]\\
=&\ch_\f\otimes\xi_\f\otimes\xi_\f\otimes\cdots
\sim\ch_\f\,.
\end{align*}
Hence, $\big[(\pi_\f(\gb)'\xi_\f\big]=\ch_\f$, and thus $\f$ has central support in the bidual.

Suppose now that $\f$ has central support in the bidual, that is $\big[(\pi_\f(\gb)'\xi_\f\big]=\ch_\f$. Then, easily
$$
\big[\pi_{\om}(\ga)'_o(\xi_\f\otimes\xi_\f\otimes\xi_\f\otimes\cdots)\big]
=\ch_{\om}\,,
$$
which means that $\om$ has central support in the bidual.
\end{proof}
%\begin{rem}
%We remark that (see \cite{vN, St}), since $\psi_{\f,\f,\dots}\cs(\ga_{\rm X})$ is a factor state if and only if $\f\cs(\gb)$ is, by applying Theorem
%\ref{symis} we conclude that $\times_\bn\f\in\cs_+(\ga_{\rm X})$ is a factor state if and only if $\f\incs_+(\gb)$ is.
%\end{rem}

To conclude, we list a class of results for the symmetric states on the Fermi algebra $\ga_{\rm F}$, corresponding to the analogous ones for symmetric states on $\ga_{\rm X}$ in \cite{St2}, which the proofs again follow by the simple application of Theorem \ref{klete}.
\begin{prop}[\cite{St2}, Thm. 2.2]
Let $\f\in\cs_+(\gb)$ be an even state. The product state $\times_\bn\f\in\cs_+(\ga_{\rm X})$ is a factor state if and only if $\f$ is. 

In addition,
\begin{itemize}
\item[(a)] $\times_\bn\f$ is of type $\ty{I}_1$ if and only if $\f$ is a character;
\item[(b)] $\times_\bn\f$ is of type $\ty{I}_\infty$ if and only if $\f$ is pure, but not  a character;
\item[(c)] $\times_\bn\f$ is of type $\ty{II}_1$ if and only if $\f$ is trace, but not  a character;
\item[(d)] $\times_\bn\f$ is of type $\ty{II}_\infty$ if and only if $\om_{\xi_\f}\lceil_{\pi_\f(\gb)'}$ is trace, and $\f$ is neither pure nor a trace;
\item[(e)] $\times_\bn\f$ is of type $\ty{III}$ if and only if $\om_{\xi_\f}\lceil_{\pi_\f(\gb)'}$ is not a trace.
\end{itemize}
\end{prop}
If $\gb\sim\bm_n(\bc)$ (we simply write $\gb=\bm_n(\bc)$) is a finite dimensional von Neumann factor, any state $\f\in\cs(\gb)$ is uniquely given by a positive, normalised "density matrix" $\r_\f\in\bm_n(\bc)$ such that $\f(a)=\tr_n(\r_\f a)$. It is even if and only if $\r_\f$ is.\footnote{It is now clearly explained why, in order to define a product state on the Fermi product of finite dimensional $C^*$-algebras (or in general of semifinite $W^*$-algebras), it is necessary and sufficient that at least one of them should be even. In fact, the density ({\it i.e.} the Radon-Nikodym derivative) of the product state w.r.t. the canonical trace should be the product of the corresponding densities, provided it commute each other. But this can happen if and only if at least one of the involved densities is even.}

In this case, the factor associated to a product state on $\ga_{\rm X}$, generated by a single state, is a particular case of an {\it Infinite Tensor Product of  
factors of finite type} $\ty{I}$ (ITPFI for short).  

%The following result covers Corollary 2.4 of \cite{St2}.
\begin{rem}[\cite{St2}, Corollary 2.4]
The product state $\times_\bn\f\in\cs(\ga_{\rm F})$, $\f\in\cs_+(\gb)$, is $*$-isomorphic to an {\rm ITPFI} factor, and thus the type of the arising factor is uniquely determined by the eigenvalue-list of the density $\r_\f$, as made explicit in Corollary 2.4 of \cite{St2}.\footnote{See also \cite{FM1}, Theorem 5.3, for a similar situation arising from the more complicated situation of the so-called quantum Markov states.}
\end{rem}

\begin{prop}[\cite{St2}, Thm. 2.11] 
Let $\gb$ be a von Neumann factor different from the scalars. Then, with the notations at the beginning of the present section, 
$\cs_\bp(\ga)_{\rm III}$ is dense in $\cs_\bp(\ga)$. If in addition, $\gb$ is not of type $\ty{I}$, then also $\cs_\bp(\ga)_{\rm I}$ is dense in $\cs_\bp(\ga)$.
\end{prop}

\section{A simple illustrative example}
\label{ezx}

We describe an example which, perhaps quite simple, helps to understand that the situation arising from product states of Fermi algebras is rather more complicated than that of the usual tensor product. This is essentially due to the strong entanglement arising from the (anti)commutation relations of the Fermi algebras. 

The first crucial difference is that, in general,
\begin{equation}
\label{nocom}
Z\big(\cam\, \overline{\circled{\rm{{\tiny F}}}}\, \cn\big)
\neq
Z(\cam)\, \overline{\circled{\rm{{\tiny F}}}}\, Z(\cn)\,,
\end{equation}
where $\cam$ and $\cn$ are two graded von Neumann algebras.

In fact, the 1st algebra is abelian by construction, whereas the 2nd one is not in general. To immediately recognise the contradiction, suppose that 
$\cam\equiv\ca$ and 
$\cn\equiv\cb$ are %Maximal Abelian Sub-Algebras (MASA for short) 
abelian, and assume that 
$Z(\ca)\, \overline{\circled{\rm{{\tiny F}}}}\, Z(\cb)=Z\big(\ca\, \overline{\circled{\rm{{\tiny F}}}}\, \cb\big)$.
We would get 
$$
\ca\, \overline{\circled{\rm{{\tiny F}}}}\,\cb\equiv Z(\ca)\, \overline{\circled{\rm{{\tiny F}}}}\, Z(\cb)
=Z\big(\ca\, \overline{\circled{\rm{{\tiny F}}}}\, \cb\big)\subset(\ca\, \overline{\circled{\rm{{\tiny F}}}}\,\cb)'\,,
$$
which means that $\ca\, \overline{\circled{\rm{{\tiny F}}}}\,\cb$ would be abelian as well. As we are going to show it is not true in general.

Indeed, consider the simplest CAR algebra $\big({\rm CAR}(\{1,2\}),\th\big)$ generated by two independent annihilators $a, A$, where $a\sim\begin{pmatrix} 
	 0&1\\
	0& 0\\
     \end{pmatrix}$, $\th$ being the usual grading automorphism of the CAR algebra. Define $s:=a+a^\dagger$, $S:=A+A^\dagger$ the corresponding "position operators", and consider the two abelian algebras 
$$
\ca:=\{\a p+\b q\mid\a,\b\in\bc\}\,,\quad  \cb:=\{\a P+\b Q\mid\a,\b\in\bc\} \,,
$$ 
generated by the projections
$$
p:=\frac12(I+s)\,,\,\,\,q:=\frac12(I-s)\,;\quad P:=\frac12(I+S)\,,\,\,\,Q:=\frac12(I-S)\,,
$$
where $I=\idd_{{\rm CAR}(\{1,2\})}$. The details of such a framework are described in \cite{AFM}, Proposition 4.4. Here, we only note that, since the above projections are not homogeneous ({\it i.e.} neither even nor odd), $\th(p)=q$, and the analogous one $\th(P)=Q$, hold true.

Therefore, $\big(\ca,\th\lceil_\ca\big)$ and $\big(\cb,\th\lceil_\cb\big)$ are $\bz^2$-graded $C^*$-algebras, where perhaps the grading is obviously not inner, and it is
straightforward to convince oneself that, inside ${\rm CAR}(\{1,2\})$, 
$$
{\rm CAR}(\{1,2\})\supset\ca\vee\cb\sim\ca\circled{\rm{{\tiny F}}}\cb\equiv\ca\circled{\rm{{\tiny F}}}\ca\,,
$$
with $Z\big(\ca)=\ca$, $Z\big(\cb)=\cb$.\footnote{Since $\dim\big({\rm CAR}(\{1,2\})\big)=16$, the generated von Neumann algebra is simply the $*$-algebra generated by $\ga$ and $\cb$ in ${\rm CAR}(\{1,2\})$.} On the other hand, 
$[p,P]=\frac{sS}2$, after an easy computation which uses the Fermi anticommutation rules.

Therefore, $Z\big(\ca\, \overline{\circled{\rm{{\tiny F}}}}\, \cb\big)$ is abelian but $Z(\ca)\, \overline{\circled{\rm{{\tiny F}}}}\, Z(\cb)=
\ca\, \overline{\circled{\rm{{\tiny F}}}}\, \cb$ is not, whence \eqref{nocom}.

We now report all the emerging structure for this simple example. 
First, it is quite straightforward to recognise that $\cam\sim\bm_2(\bc)$. In fact, by using the commutation relations, 
$\cam={\rm span}\big\{I,s,S,sS\big\}$, and thus it is a 4-dimensional non abelian $C^*$-algebra, whence the claim.

It is also quite easy to see that
\begin{equation}
\label{matfrs}
e_{11}:=p\,,\,\,e_{22}:=q\,,\,\,e_{12}:=pS\,,\,\,e_{21}:=qS
\end{equation}
provides a system of matrix-units for $\cam$.

Concerning the even product states $\f$, there is only one on $\ca$: $\f(p)=1/2=\f(q)$. In the matrix-representation using \eqref{matfrs},
it is easy to see that the product state $\om=\f\times\f$ is given by 
$$
\om\Big(\sum_{i,j=1}^2 a_{ij}e_{ij}\Big)=\frac12\tr\begin{pmatrix} 
	 a_{11}&a_{12}\\
	a_{21}&a_{22}\\
     \end{pmatrix}\,.
$$
Therefore,
\begin{align*}
\bm_2(\bc)\sim\pi_{\f\times\f}(\ca\circled{\rm{{\tiny F}}}\ca)\equiv&\pi_{\f\times\f}(\ca\,\circled{\rm{{\tiny F}}}\,\ca)''
\nsim\pi_{\psi_{\f,\f}}(\ca\circled{\rm{{\tiny X}}}\ca)''\\
\equiv&\pi_{\psi_{\f,\f}}(\ca\circled{\rm{{\tiny X}}}\ca)\sim\bc^2\otimes\bc^2
\end{align*}
and, in particular, the factor representation $\pi_{\f\times\f}$ of the product state $\f\times\f$ is generated by the non factor representation of $\pi_\f$ of the single algebra $\ca$.

To comment the general case, we consider two arbitrary $\bz^2$-graded $C^*$-algebras $\ga$ and $\gb$, together with two even states $\om$ and $\f$ on $\ga$ and $\gb$, respectively. 

The example just explained is affirming that, in general, there is no Klein transformation mapping $\ga\,\circled{\rm{{\tiny F}}}\,\gb$ isomorphically onto the corresponding object $\ga\,\circled{\rm{{\tiny X}}}\,\gb$, or even for fixed representations $\pi_\om(\ga)''\,\circled{\rm{{\tiny F}}}\,\pi_\f(\gb)''$
and the corresponding ones $\pi_\om(\ga)''\,\circled{\rm{{\tiny F}}}\,\pi_\f(\gb)''$. This simply happens because the Fermi and the usual tensor product generate, in general, not isomorphic $C^*$-algebras as we have seen above.

In particular, using the notations in Section 5 of \cite{CDF}, the computation 
\begin{align*}
{\rm Z}_{\om\times\f}\sim&\eta_{V_\om\otimes V_\f}\big(\overline{\pi_{\om\times\f}(\ga\,\circled{\rm{{\tiny F}}}\,\gb)}\big)
\bigcap\big(\overline{\pi_{\om\times\f}(\ga\,\circled{\rm{{\tiny F}}}\,\gb)}\big)^\wr\\
=&\Big(\eta_{V_\om}\big(\overline{\pi_{\om}(\ga)}\big)
\bigcap\big(\overline{\pi_{\om}(\ga)}\big)^\wr\Big)\,\circled{\rm{{\tiny F}}}\,\Big(\eta_{V_\f}\big(\overline{\pi_{\f}(\ga)}\big)
\bigcap\big(\overline{\pi_{\f}(\ga)}\big)^\wr\Big)\\
\sim&{\rm Z}_\om\,\overline{\circled{\rm{{\tiny F}}}}\,{\rm Z}_\f
\end{align*}
appeared in \cite{CRZ2} is incorrect.\footnote{The above example shows also that, in general, ${\rm Z}_{\om\times\f}\nsim{\rm Z}_\om\,\overline{\circled{\rm{{\tiny X}}}}\,{\rm Z}_\f\sim
{\rm Z}_{\psi_{\om,\f}}$}

To conclude, even if the simplex $\cs_\bp(\ga_{\rm F})$ of the symmetric states on the Fermi algebra $\ga_{\rm F}$ shares many properties with the corresponding object $\cs_{\bp\times\bz^2}(\ga_{\rm X})$ of the even and symmetric states on the usual tensor product algebra $\ga_{\rm X}$, the properties of the corresponding elements ({\it i.e.} Fermi product states on $\ga_{\rm F}$ and the corresponding product states on $\ga_{\rm X}$) can be never compared in general.

\end{document}